\documentclass[a4paper,11pt,reqno]{amsart}

\usepackage{amssymb,amsthm,amsmath}
\newtheorem{secc}{}
                         \newtheorem{theorem}{Theorem}

                         \newtheorem{prop}{Proposition}

                         \newtheorem{lemma}{Lemma}

                         \newtheorem{rema}{Remark}

                         \newtheorem{defn}{Definition}

                         \newtheorem{cor}{Corollary}

\def\a{\theta}
\newcommand{\T}{\mathbb{T}}
\newcommand{\C}{\mathbb{C}}
\newcommand{\R}{\mathbb{R}}
\newcommand{\Q}{\mathbb{Q}}
\newcommand{\Z}{\mathbb{Z}}
\newcommand{\N}{\mathbb{N}}

\newcommand{\eps}{\varepsilon}

\def\carre{ \hfill $\Box$    }
\begin{document}
\title{On the ergodicity of the Weyl sums cocycle.}
\author{Bassam Fayad}
\address{Bassam  Fayad, LAGA, Universit\'e Paris 13, Villetaneuse, }
\email{fayadb@math.univ-paris13.fr}


\date{March 3, 2005}

\maketitle

\begin{abstract} For $\theta \in [0,1]$, we consider the map $T_\a: \T^2 \to \T^2$ given by $T_\theta(x,y)=(x+\theta,y+2x+\theta)$. The skew product $f_\a: \T^2 \times \C \to \T^2 \times \C$ given by  $f_\theta(x,y,z)=(T_\theta(x,y),z+e^{2 \pi i y})$ generates the so called Weyl sums cocycle $a_\a(x,n) = \sum_{k=0}^{n-1} e^{2\pi i(k^2\theta+kx)}$ since the $n^{{\rm th}}$ iterate of $f_\a$ writes as $f_\a^n(x,y,z)=(T_\a^n(x,y),z+e^{2\pi iy} a_\a(2x,n))$.

In this note, we improve the study developed by Forrest in  \cite{forrest2,forrest}  around the density for $x \in \T$ of the complex sequence ${\{a_\a(x,n)\}}_{n\in \N}$, by proving the ergodicity of $f_\theta$ for a class of numbers $\a$ that contains a residual set of positive Hausdorff dimension in $[0,1]$. The ergodicity of $f_\a$ implies the existence of a residual set of full Haar measure of $x \in \T$ for which the sequence ${\{ a_\a(x,n) \}}_{n \in \N}$ is dense.
\end{abstract}

\begin{secc}\rm Let $\T^2$ denote the torus $\R^2 / \Z^2$. For $\theta \in [0,1]$ define the map ({\it skew shift}) $T_\theta$:
\begin{eqnarray*} \T^2 &\to& \T^2 \\
(x,y) &\mapsto& (x+\theta,y+2x+\theta) \end{eqnarray*}
and the {\it skew product}  $f_\a$:
\begin{eqnarray*} \T^2 \times \C &\to& \T^2 \times \C \\
(x,y,z) &\mapsto& (x+\theta,y+2x+\theta,z+e(y)) \end{eqnarray*}
where $e(y)$ is the usual notation for $e^{2 \pi i y}$. The diffeomorphism $f_\theta$ preserves the product measure $\mu=m \times \nu$ where $m$ denotes the Haar measure on $\T^2$ and $\nu$ denotes the Lebesgue measure on $\C$. We say that the map $f_\a$ is {\it ergodic} if and only if for every $\mu$-measurable set $A \subset \T^2\times \C$ such that $f_\a(A)=A$ we have $\mu(A)=0$ or $\mu(A^c)=0$.

\begin{defn}\rm We define ${\mathcal F}$ to be the set of numbers  $\theta \in [0,1] \setminus \Q$ having a continued fraction representation 
$$
\a=\frac{1}{\displaystyle{a_1+\frac{1}{a_2+\frac1{\ldots}}}},
$$

\vspace{0.3cm}

\noindent such that $\sum_n 1/a_n < \infty$, and such that $\liminf_{q\geq 1} q^{3 +\eps} \|q\theta\|=0$ for some $\eps>0$. Here and in all the text $\| \ldots \|$ stands for the closest distance of a real number to the integers. Let $p_l/q_l = [a_1,\ldots,a_l] = 1/(a_1+1/(a_2+\ldots+(1+ 1/a_l)\ldots))$, $p_l$ and $q_l$ relatively prime. The sequence $p_l/q_l$ is called the sequence of the best rational approximations of $\a$ since we have $\|q_{l-1}\a\| \leq q_{l-1}\a$ for every $k < q_{l}$. The sequence $q_l$ is simply called the sequence of {\it approximation denominators} of $\a$.
\end{defn}

We will elaborate on the paper by Forrest \cite{forrest} to obtain the following result: 
\begin{theorem} \label{main} Let $\a \in {\mathcal F}$. Then $f_\theta$ is ergodic.
\end{theorem}

The set ${\mathcal F}$ has zero measure due to any of the two conditions imposed on $\a$. 
The set has positive Hausdorff dimension but the condition $\sum 1/a_n < \infty$ on $\a$ is actually very restrictive since it involves {\it all} the convergents of $\a$. For instance ${\mathcal F}$ is contained in the complementary of a residual set (this can be checked by the ergodicity of the Gauss transformation $\theta \mapsto \{1 / \theta \}$). But we can show using  a classical general argument of Halmos, exposed in his introductory book to ergodic theory \cite[proof of the second category theorem]{halmos}, that the set of $\a$ such that $f_\a$ is ergodic is a $G_\delta$ set, call it  $\tilde{{\mathcal F}}$. Since ${\mathcal F}$ is dense and ${\mathcal F} \subset \tilde{{\mathcal F}}$, we have

\begin{cor} \label{22} The set $\tilde{{\mathcal F}} \subset [0,1]$ of $\a$ such that $f_\a$ is ergodic is a residual set of positive Hausdorff dimension.
\end{cor}

This actually hints at the possibility of bypassing the condition $\sum 1/a_n < \infty$ in the proof of ergodicity.  Proposition \ref{pop} and hence proposition \ref{for}, that are the only places where this condition appears, can actually be proven without it using recent results on theta sums. This will be done in a future work. 

\end{secc}

\begin{secc}\rm Theorem \ref{main} and its corollary are a strengthening of the main result of \cite{forrest} where the density in $\C$ of the Weyl sums 
  \begin{eqnarray} \label{wsum} \sum_{k=0}^{n-1} e(k^2\theta+kx),  \quad n=1,2,\ldots \end{eqnarray}
was proved, if $\a \in {\mathcal F}$, for almost every $x \in [0,1]$. Indeed, we have
\begin{cor} Let $\a \in \tilde{{\mathcal F}}$. Then the set
$$B(\theta) = \lbrace x \in [0,1] : \sum_{k=0}^{n-1} e(k^2\theta+kx),  \  n=1,2,\ldots, \text{is dense in} \ \C \rbrace $$
is a $G_\delta$ dense set of full Lebesgue measure in $[0,1]$.  
\end{cor}
\begin{proof} If $f_\a$ is ergodic then for $\mu$-a.e. $u=(x,y,z)$ we have that the sequence $u,f_\a(u), f_\a^2(u),\ldots,$ is dense in $\T^2 \times \C$. This is indeed a general fact that can be proved considering a countable base ${\{O_j\}}_{ j \in \N}$ of open balls of $\T^2 \times \C$ and observing that the complementary of the invariant set $\cup_{n \in \Z} f_\a^n(O_j)$ has zero measure  from which it follows that the complemetary of the set ${\mathcal D}= \cap_{j \in \N} \cup_{n \in \Z} f_\a^n(O_j)$. But by definition a point $x \in {\mathcal D}$ has a dense orbit under $f_\a$.   Now 
\begin{eqnarray} \label{iterates} f_\a^n(x,y,z)=(T_\a^n(x,y), z+\sum_{k=0}^{n-1}e(k^2\a+2kx+y)), \end{eqnarray}
so that for $\mu$-a.e. $(x,y,z)$ we have that the sequence  
 $ z+\sum_{k=0}^{n-1}e(k^2\a+2kx+y), \  n=1,2,\ldots, $ is dense in $\C $. The density of the latter sequence clearly does not depend on $y$ and $z$ and the measurable statement of the corollary follows. Further, ${\mathcal D}$ is a $G_\delta$ set and since its complementary has zero measure it follows that it is a dense $G_\delta$. For the same reason as above this means that $B(\theta)$ is a $G_\delta$-dense set. \end{proof}
 
We will see that in  proving the density of the Weyl sums (\ref{wsum}) for almost every $x$ when $\theta \in {\mathcal F}$, Forrest actually went a long way towards proving the ergodicity of $f_\a$. Yet, he left this question unsolved and put it as an open problem even for a single value of $\theta$. In a sense, we will finish here his work.

Finally, we recall that prior to \cite{forrest}, Forrest had already proved in \cite{forrest2} the transitivity of $f_\a$ under the sole hypothesis $\liminf_{q \geq 1} q^{3/2} \|q \a\| <\infty$. From the transitivity of $f_\a$, the density of the Weyl sums follows for
 a dense $G_\delta$ set of $x\in[0,1]$. Although $T_\a$ is uniquely ergodic, the cocycles $ \sum_{k=0}^{n-1}e(k^2\a+2kx+y))$  behave differently for different points $(x,y) \in \T^2$ as shown by the following remark:
 
\begin{rema}\rm While it is not clear whether $0$ could be in $B(\theta)$ for some choice of $\a$\footnote{The claim made by Forrest that it follows from \cite{HL} that  $0 \notin B(\theta)$    for any irrational $\a$  probably stems from his  missinterpretation of the formula $a(0,n) = \Omega(\sqrt{n})$ which is used in \cite{HL} (cf. \S \ref{four} below) as  the negation of $a(0,n) = o(\sqrt{n})$ and not as $\sqrt{n} = O(|a(0,n)|)$ like Forrest might have understood it. It is clear from the formulae of $a(0,n)$ in the case of $\a$ rational that one can construct an irrational $\a$ for which there exists a seqeunce  $q_n \to \infty$ such that $a(0,q_n) \to 0$.}, it does follow from 
an argument by Besicovitch \cite{besicovitch} that for any $\theta$ there exists always an $x$ such that $x \notin B(\theta)$.
\end{rema}
\end{secc}

\begin{secc}\rm The question of knowing whether the set ${\mathcal F}$ of $\a$ for which $f_\a$ is ergodic (or even transitive) has full measure (or contains all irrationals!) is still open
 and we have not much to say about this as explained in the following list of remarks:

\begin{rema}\rm  \label{rem}
It does not seem to be known whether there exists a class of irrational numbers $\a$  for which the Weyl sums could fail to be dense for every $x$. In \cite{forrest} it is claimed erroneously\footnote{For the same reason as in the precedent footnote.} that  the estimate $|\sum_{k=0}^{n-1} e(k^2\theta+kx)  | \geq c_\a \sqrt{n}$ (uniformly in $x \in [0,1]$) was proved in \cite{HL} for constant type numbers $\a$ (numbers with bounded partial quotients, or equivalently numbers that satisfy $\liminf_{q \geq 1} q \|q\a\| >0$). If this however turns out to be true, it would obviously preclude, if $\a$ is of constant type, the density of the Weyl sums for any choice of $x$.

Remarkably, if true, the latter estimate turns out to be paradoxically helpful in showing ergodicity of the Weyl sums without the restrictive hypothesis $\sum 1/a_n<\infty$. Indeed,
 an elegant proof of ergodicity of $f_\a$ for some class of $\a$ (included in those satisfying $\liminf_{q \geq 1} q^5 \|q \a\|=0$)  was given in \cite{volny}, that is based on the alleged uniform lower bound on the Weyl sums for constant type numbers $\a$.
\end{rema}

\begin{rema}\rm While a property on the rational approximations of $\a$, at least like the one used in \cite{forrest2}, namely $\lim \inf q^{3/2}\|q\a\| = 0$, seems necessary to  study the density of the Weyl sums using the dynamics of $f_\a$, the condition $\sum_{n \geq 1}  1/a_n < + \infty$ could be removed as in \cite{volny} from the proof if some upper bounds on the measure of the sets where $|\sum_{k=0}^{n-1} e(k^2\theta+kx)|$ are not large  were known. It would be helpful for example if one knows that for any constant $C>0$,
$$\lim_{q \to \infty} \sup_{1 \leq p \leq q-1} \lambda \{ x : |\sum_{k=0}^{q-1} e(k^2 p/q + kx)| \leq C \}=0.$$  \end{rema}
 
\begin{rema}\rm If we denote for $l \geq 1$ by $f_\a^{(l)}$ the skew product $f_\theta^{(l)}(x,y,z)=(x+\theta,y+2x+\theta,z+e(ly))$, then the same proof of ergodicity for $\a \in {\mathcal F}$ of $f_\a^{(1)}$ implies the ergodicity of every $f_\a^{(l)}$. But the set of $\a \in [0,1]$ with the latter property is invariant by multiplication by $l$ on the circle so has measure either $0$ or $1$. 

To compare with our problem, note that twist maps of the type $\T^d \times \R^k \to \T^d \times \R^k$, $(x,z) \to (x+\alpha,z+\varphi(x))$ with a smooth function $\varphi$ having zero average and that is not a trigonometric polynomial are always ergodic for a $G_\delta$-dense set of $\alpha \in \T^d$ (of zero Hausdorff dimension however) and not ergodic for a set of $\alpha$ of full measure which  consists of the Diophantine vectors, that is vectors for which there exists $N$ such that $\liminf_{q \geq 1} q^N \|q\alpha \| > 0$.   
\end{rema}
\end{secc}

\begin{secc}\rm \label{four} In \cite{HL}, Hardy and Littlewood studied the growth of $|\sum_{k=0}^{n-1} e(k^2\theta+kx)|$ for different values of $\a \in [0,1]$. Using the notation $u_n=\Omega(v_n)$ for positive sequences $u_n$ and $v_n$ for the negation of $u_n=o(v_n)$,
 the principal bounds they obtained were 
 
\vspace{0.2cm}

\noindent {\bf Theorem} \cite[Theorems 2.14, 2.141, 2.18, 2.181, 2.22, 2.221]{HL} {\sl For any irrational $\a \in [0,1]$, 
$$\left|\sum_{k=0}^{n-1} e(k^2\theta+kx)\right| =o(n), \quad {\rm uniformly \ for \ all \ values \ of \ } x.$$ 
If the partial quotients $a_n$ in the continued fraction expansion of $\a$ are bounded then
$$\left|\sum_{k=0}^{n-1} e(k^2\theta+kx)\right| =O(\sqrt{n}), \quad {\rm uniformly \ for \ all \ values \ of \ } x.$$
These are optimal bounds. Indeed, for any irrational $\a \in [0,1]$ we have
$$\left|\sum_{k=0}^{n-1} e(k^2\a)\right|=\Omega(\sqrt{n}),$$ 
and for every sequence $\varphi_n>0$ tending to $0$ as $n \to \infty$, it is possible to find irrationals $\a$ such that 
$$\left|\sum_{k=0}^{n-1} e(k^2\a)\right|=\Omega(n\varphi_n).$$}

\vspace{0.2cm}

With the dynamical approach adopted in this paper, the first one of these equations follows immediately from two classical and elementary facts in ergodic theory, see e.g. \cite{parry}: first, that $T_\a$ is uniquely ergodic as soon as $\a$ is irrational; and second, that this implies that the function $\Phi(x,y)= e(y)$, of zero average, has its Birkhoff means $1/n \sum_{k=0}^{n-1} e(k^2\theta+2kx+y)$ converging uniformly to zero.

It would be nice if a an additional qualitative ergodic property of $T_\a$ could be  displayed in the case of irrationals $\a$ with bounded partial quotient that would explain the second bound in the above theorem of Hardy and Littlewood. 

\end{secc}

\begin{secc}\rm We now procede to the proof of theorem \ref{main}. In all the sequel, $\a$ will be a fixed irrational number in ${\mathcal F}$. For every $n,m\in \N$ and $(x,y) \in \T^2$, let 
$$a(x,y,n)= \sum_{k=0}^{n-1} e(k^2\theta + 2kx+y),$$
and  
$$b(x,m)= \sum_{k=0}^{m-1} e(kx).$$
 
\begin{defn}\rm[Essential value] 
We say that $l \in \C$ is an {\it essential value} for the cocycle $a$ above $T_\a$ if for any measurable set $E \subset \T^2$ such that $m(E)>0$ and for any $\nu>0$, there exists $n \in \N$ such that 
$$m\left( E \cap T_\theta^{-n} E \cap \lbrace (x,y) \ / \ \left| a(x,y,n)-l \right| \leq \nu \rbrace  \right)>0.$$
We say that $l\geq 0$ is an {\it essential value for the modulus of $a$}  if 
for any measurable set $E \subset \T^2$ such that $m(E)>0$ and for any $\nu>0$, there exists $n \in \N$ such that 
$$m\left( E \cap T_\theta^{-n} E \cap \lbrace (x,y) \ / \ \left| |a(x,y,n)|-l \right| \leq \nu \rbrace  \right)>0.$$   
Since $|a(x,y,n)|$ does not depend on $y$ we simply denote it by $|a(x,n)|$.
\end{defn}

A very useful general criterion for ergodicity established by K. Schmidt in \cite{schmidt} states that $f_\a$ is ergodic if and only if any $l \in \C$ is an essential value for $a$ (above $T_\a$), but due to the symmetries of the system we have the following sufficient criterion for ergodicity that we took from \cite{volny}:
\begin{lemma} \label{mod}
If $1/2$ (or any other strictly positive number) is an essential value for the modulus of $a$ then $f_\a$ is ergodic.
\end{lemma}

\begin{proof} The proof contains two parts. First, it is shown that $a$ has a nonzero 
essential value. Indeed, if this was not true, then by lemma \cite[Lemma 3.8]{schmidt} (the proof of this lemma can also be found in \cite[Lemma 8.4.3]{aaronson}), we have that for any compact set $K \subset \C$ that does not contain $0$, there exists a measurable set $B \subset \T^2$ such that for every $n \in \N$,
$$B \cap T_\theta^{-n} B \cap \lbrace (x,y) \ / \ a(x,y,n) \in K  \rbrace = \emptyset$$
which clearly contradicts the assumption of the lemma.

Next, assume that $l \neq 0$ is an essential value for $a$. For $y_0 \in \T$ denote by $S_{y_0}$ the map of $\T^2$ on itself $S_{y_0} (x,y)=(x,y+y_0)$. Then, the fact that for a measurable set $B$ with $m(B) >0$ we have an $n \in \N$ such that
$$m\left( S_{y_0}B \cap T_\theta^{-n} (S_{y_0}B) \cap \lbrace (x,y) \ / \ \left| a(x,y,n)-l \right| \leq \nu \rbrace  \right)>0,$$   
 implies for the same $n$ that 
$$m\left(B \cap T_\theta^{-n} B \cap \lbrace (x,y) \ / \ \left| a(x,y,n)-l e(-y_0) \right| \leq \nu \rbrace  \right)>0,$$  
which implies that all the circle of modulus $|l|$ is included in the set of essentail values of $a$. Since the set of essential values of a complex cocycle above an ergodic map is a closed subgroup of $\C$ (cf. \cite[Lemma 3.3]{schmidt}), it follows that for the cocycle $a$ it is equal to $\C$ and $f_\a$ is hence ergodic. 
\end{proof}
\end{secc}

\begin{secc}\rm The general strategy in controlling $|a(x,n)| = |\sum_{k=0}^{n-1} e(k^2\theta + 2kx)|$  starts by showing that given any infinite subsequence of the approximation denominators of $\a$, and in particular along a subsequence that satisfies the hypothesis  $q_n^{3+\eps}\|q_n\a\| \to 0$, we have that  for a typical value of $x$, $|a(x,q_n)| \to \infty$. This implies an approximation formula for $|a(x,mq_n)|$, when $m$ is not too large, by $|a(x,q_n)||b(2q_n x,m)|$  and $m$ is then chosen to bring this product close  to $1/2$. Typically, when $2q_nx$ behaves like a badly approximated number, $|b(2q_nx,l)|$, $l=1, \ldots, m$  contains a $O(1/m^{1-\epsilon})$-dense set in $[0,1]$ (here $\epsilon >0$ is an arbitrarilly small number). If we prove that $|a(x,q_n)|$ is typically bounded by $q_n^{1/2+\epsilon}$ then the $m_n$ we need to modulate the product $|a(x,q_n)||b(2q_n x,m)|$  is not larger than  $q_n^{1/2+2\epsilon}$ and the condition $q_n^{3+\eps}\|q_n\a\| \to 0$ appears then to be the exact condition that allows the approximation formula to hold up to this value of $m$. 

Finally, to show that $1/2$ is actually an essential value for the modulus of $a$ we compute a bound on the derivative with respect to $x$ of the product $|a(x,q_n)||b(2q_n x,m_n)|$ and show that, under the  same assumption $q_n^{3+\eps}\|q_n\a\| \to 0$, the interval $I_n$ containing $x$ where the product is close to $1/2$ is sufficiently large so that $R_\a^{m_nq_n}(I_n)$  is almost equal to $I_n$. This and the fact that $|a(x,y,l)|$ does not depend on $y$ will  allow us to conclude.

In this scheme, the first step is the most delicate. It was proved by Forrest in \cite{forrest} who based his proof on the following approximate functional equation, established by Hardy and Littlewood in \cite[Theorem 2.128, Theorem 2.17]{HL}:  for $0<\a$, $x<1$ and $k\geq 1$ 
\begin{eqnarray} \label{oo} \sqrt{\a} |a(\a/2,x/2,k)| = |a(\{{1 / \theta} \}/2, \{{-x / \a }\}/2, [k\a])| +O(1) \end{eqnarray}
where $\{\cdot\}$ and $[\cdot]$ denote the fractional and the integer part of a number and where the constant involved in the Landau's error notation is absolute. Under  an additional assumption on $\a$ it is possible to apply a dynamical approach  where $\a$ is viewed as a parameter and obtain by induction from the above functional equation a lower estimate on the Weyl sums. The upshot of this approach is the following key ingredient of \cite{forrest} as well as for us here:
\begin{prop}\cite[Proposition 4.3]{forrest} \label{for}
 Suppose $\theta \in [0,1] \setminus \Q$ has a continued fraction representation $[a_1,a_2,\ldots]$ such that $\sum_n 1/a_n < \infty$. Then, given any $\delta>0$ and any infinite subset $Q$ of the set of approximation denominators of $\a$  we have that for Lebesgue almost  every $x \in [0,1]$, there exists a sequence $q_n \in Q$ such that $\delta/2 \leq \|2q_nx\|\leq \delta$  and $\mathop{\lim} \limits_{n \to \infty} |a(x,q_n)|=\infty$. 
\end{prop}

For the commoditiy of the reader and to keep this note as much self contained as possible (modulo the functional equation (\ref{oo}) that is admitted), we include in an appendix the scheme of the proof given in \cite{forrest} of the above proposition.

\end{secc}

\begin{secc}\rm To proceed we need the following construction similar to the one made in \cite{forrest}. Suppose $\a \in {\mathcal F}$, then there exists a sequence $q_n$ of approximation denominators of $\a$ such that:

\vspace{0.1cm}

 \noindent {\bf 7.a.} $q_n^{3+\eps}\|q_n\a\| \longrightarrow 0$.

\vspace{0.2cm}

\noindent {\bf 7.b.} For almost every $x \in [0,1]$ there is a sequence $U_n \to \infty$ and infinitely many $n$ such that $\delta/2 \leq \|2q_nx\|\leq \delta$ and $|a(x,q_n)| \geq U_n$ (this is exactly proposition \ref{for}).  

\vspace{0.2cm}

\noindent {\bf 7.c.} For almost every $x\in [0,1]$, there is an $n_1$ such that for $n\geq n_1$, we have $|a(x,q_n)| \leq q_n^{1/2 + \eps/10}$.

\vspace{0.1cm} 
This is because of the fact that $\int_0^1 {|a(x,q_n)|}^2 dx = q_n$ implies $\lambda \lbrace x : |a(x,q_n)| \geq q_n^{1/2+\eps/10} \rbrace \leq 1/q_n^{\eps/5}$; but 7.a implies that $q_{n+1} \geq q_n^3$, hence $ \sum 1/q_n^{\eps/5} < \infty$ and $3$ follows by the Borel Cantelli lemma.

\vspace{0.2cm}

\noindent {\bf 7.d.} For almost every $x \in [0,1]$, there is an $n_2$ such that for $n \geq n_2$, the set $\lbrace |b(2q_nx,m)| : 0 \leq m \leq q_n^{1/2+\eps/4} \rbrace$ is $1/(q_n^{1/2+\eps/8} \|2q_nx\|)$-dense in $[0,1]$.
 
\vspace{0.1cm} 

To prove this we define $H_n:=q_n^{1/2+\eps/4}$. We let $A_k^\eps \subset [0,1]$ be the subset of irrationals such that for each $\alpha \in A_k^\eps$, and for $m\geq k$, there exists a continued fraction approximation $p/q$ for $\alpha$ such that $q \in [m^{1-\eps/10},m]$. Since the set of numbers $\alpha$ for which there exists  $C>0$ such that $q_{n+1}(\alpha)  \leq {q_n(\alpha)}^{1+\varepsilon/10}$ is of full measure, we clearly have $\lambda(\cup_k A_k^\eps)=1$ and we pose $\lambda(A_k^\eps) = 1-\upsilon(k)$. In our choice of the sequence $q_n$ in  7.a we can assume up to extracting that $\sum_n \upsilon({H_n}) < \infty$. Since $\lambda \lbrace x \in [0,1] : 2q_nx {\rm  \ mod \ } [1] \in A_{H_n}^\eps \rbrace = \lambda (A_{H_n}^\eps) = 1 - \upsilon({H_n})$ we deduce that for almost every $x \in [0,1]$, there exists $n_2$ such that for $n\geq n_2$, then  $2q_nx {\rm \ mod \ } [1] \in A_{H_n}^\eps$ from which 7.d follows easily.

\vspace{0.2cm}

\end{secc}

\begin{secc}\rm Note that a simple computation (see \cite[Lemma A.4]{forrest}) gives that for some constant $C$ and for any $x \in [0,1]$, $l,m \in \N$, we have  
$$ |a(x,ml) -a(x,l)b(2lx,m)  | \leq C|a(x,l)|m^3l\|l\a\|,$$
which  in the case of $q_n$ satisfying 7.a and $m \leq q_n^{1/2+\eps/4}$ yields
\begin{eqnarray*}
|a(x,mq_n) -a(x,q_n)b(2q_nx,m)  | \leq C|a(x,q_n)|q_n^{-1/2-\eps/4},   \end{eqnarray*}
and finally, if in addition $|a(x,q_n)| \leq 2q_n^{1/2 + \eps/10}$, then 
\begin{eqnarray} \label{approx} 
|a(x,mq_n) -a(x,q_n)b(2q_nx,m)  |  \leq C q_n^{-\eps/8}.  \end{eqnarray}
It is in the above equations that the restrictive assumption $\liminf q^{3 +\varepsilon} \|q\a\| =0$ is really crucial. 

On the other hand, we have $b(2q_nx,m)= e^{i2\pi (m-1)q_nx} \sin(2\pi m q_nx) / \sin(2 \pi q_nx)$. Hence for $\delta/4 \leq \|2q_nx \| \leq 2 \delta$  we have $|b(2q_nx,m)| \leq 1/\delta$ and $|D_x (b(2q_nx,m))| \leq 4\pi mq_n / \delta$, where $D_x$ denotes derivation with respect to $x$. Also, we clearly have $|a(x,q_n)| \leq q_n$ and $|D_x (a(x,q_n) | \leq 2\pi q_n^2$. From these observations we conclude that for $n$ sufficently large, for any $m \leq q_n^{1/2 + \eps/4}$ and $\delta/4 \leq \|2q_nx \| \leq 2 \delta$, we have
\begin{eqnarray} \label{derivee}
|D_x \left[a(x,q_n)b(2q_nx,m)\right]  | \leq {5 \pi \over \delta}  q_n^{2+1/2+ \eps/4}.
\end{eqnarray}

We deduce from 7.a to 7.d the following:
\begin{prop} \label{resume} Let $\a \in {\mathcal F}$. For almost every $x \in [0,1]$ there exists an infinite sequence of integers $M_n$ and a sequence $\epsilon_n \to 0$ such that 
\begin{itemize}
\item[(i)] $\|M_n\theta\| \leq q_n^{-(2+{1 / 2} + 3 \eps / 4)}$;
\item[(ii)] For every  $\tilde{x} \in [x- q_n^{-(2+ {1 / 2} + {\eps / 2})},x+ q_n^{-(2+ 1/2 + \eps/2)}]$, we have $ ||a(\tilde{x}, M_n)|-1/2| \leq \epsilon_n$;
\item[(iii)]  $\|M_n^2\a +2M_nx\| \leq \epsilon_n$;
\end{itemize}
\end{prop}
\begin{proof} Take a sequence $q_n$ satisfying 7.a. Take an $x$ that satisfies 7.b, 7.c and 7.d. Up to extracting from $q_n$ we have that $\delta/2 \leq \|2q_nx\|\leq \delta$  and $|a(x,q_n)| \rightarrow \infty$.  From 7.c and 7.d, we find $m_n \leq q_n^{1/2+\eps/4} $ such that $|a(x,q_n)b(2q_nx,m_n)|  \rightarrow 1/2$. Since the conditions of (\ref{approx}) are satisfied by $x$ and $m_n$, (ii) follows for the particular value $\tilde{x}=x$ if we take $M_n:=m_nq_n$. 

For $|\tilde{x} -x| \leq q_n^{-(2+ {1 / 2} + {\eps / 2})}$ we have that $\delta/4 \leq \|2q_n \tilde{x} \|\leq 2\delta$, and  since $|D_x(a(\tilde{x},q_n))| \leq 2 \pi q_n^2$, we  have from 7.c that $|a(\tilde{x},q_n)| \leq 2 q_n^{1/2 + \eps/10}$, hence (\ref{approx}) holds for $\tilde{x}$ and  for the same $m_n$ considered above.  At last, (ii) then follows from  (\ref{derivee}).

From 7.a we get (i) and the fact that $\|M_n^2 \a\| \to 0$. Finally the combination of $|a(x,q_n)| \rightarrow \infty$ and  $|a(x,q_n)b(2q_nx,m_n)| \rightarrow 1/2$ forces $|b(2q_nx,m_n)| \to 0$, hence $\|2M_nx\| =\|2m_nq_nx\| \to 0$ and (iii) is proved.
\end{proof}

\begin{rema}\rm It would be possible to insure that $|a(x,q_n)b(2q_nx,m_n)|$  stays close to $1/2$ on larger intervals than in (ii) which would allow to relax the requirement (i) and from there relax the arithmetic condition 7.a on $\a$. But this condition, as we saw, is optimal if we want to insure (\ref{approx}) without which the product  $|a(x,q_n)b(2q_nx,m_n)|$ stops being interesting to our end.
\end{rema}

\end{secc}
 
\begin{secc}{\it Proof of theorem \ref{main}.} \rm From lemma  \ref{mod} it is enough to prove that 
$1/2$ is an essential value for the modulus of $a$. 

We will use $\lambda$ and $m$ to denote respectively the Haar measure on the tori $\T^1$ and $\T^2$. Fix $E \subset \T^2$ such that $m(E)>0$.  Fix then a square $A=I \times J = [x_1,x_2] \times [y_1,y_2]$, $|x_2-x_1|=|y_2-y_1|=l>0,$ such that  $m(E \cap A) \geq {(9 / 10)} m(A)$. 
We denote $m_{E\cap A}$ the induced measure: $m_{E\cap A}(B)= m(E \cap A \cap B)$ for any Lebesgue measurable set $B \subset \T^2$. We denote by $\pi^* m_{E \cap A}$ the projected measure  given by 
$\pi^* m_{E \cap A}(K)=m_{E\cap A}(K\times \T^1)$ for any Lebesgue measurable set $K \subset \T^1$. Clearly, $\pi^* m_{E \cap A} \leq \lambda$, while $\pi^* m_{E \cap A} (I) \geq  {(9/10)} l^2$, and $\lambda(I)=l$. Hence, considering the Radon-Nikodym derivative of   $\pi_* m_{E \cap A} $ with respect to $\lambda$, we find that there exists $r_0>0$ and a set $\tilde{I} \subset I$ with $\lambda(\tilde{I})>0$, such that for any $r\leq r_0$ and for any $x \in \tilde{I}$ we have $\pi_* m_{E \cap A} ([x-r,x+r]) \geq {(4 / 5)} l 2r$, that is
$$m(\Delta(x,r) \cap E) \geq  {4\over 5} m(\Delta(x,r)),$$
where $\Delta(x,r)  = [x-r,x+r] \times J.$ 

Since $\lambda(\tilde{I})>0$, it is possible to take $x_0 \in \tilde{I}$ 
for which the statement of proposition \ref{resume} holds. Recall that for $({x},{y}) \in \T^2$ and $p \in \N$ we write indifferently $|a({x},{y},p)|$ or $|a({x},p)|$ since the modulus of $a$ does not depend on $y$. 
If we denote  $\Delta_n=  \Delta(x_0,q_n^{-(2+1/2+\eps/2)})$, we have from (ii) of proposition \ref{resume}  that 
\begin{eqnarray} \label{11} ||a({x}, {y}, M_n)|-1/2| \leq \epsilon_n, \quad  {\rm for \ all \ }  ({x},{y}) \in  \Delta_n. \end{eqnarray}
From the definition of $\tilde{I}$ we have for $n$ sufficiently large 
\begin{eqnarray} \label{33} m(\Delta_n \cap E) \geq  {4\over 5} m(\Delta_n). \end{eqnarray}

On the other hand, (i) and (iii) imply that 
\begin{eqnarray} \label{222} \lim_{n\to \infty} {m \left( T_\a^{-M_n}  \Delta_n \triangle  \Delta_n \right) \over m(\Delta_n)}=0, \end{eqnarray}
where $A \triangle B$ stands for the symmetric difference between $A$ and $B$. 

It immediately follows from (\ref{33}) and (\ref{222}) (because $4/5>1/2$) that  for $n$ sufficiently large
$$m( T_\theta^{-M_n}E \cap E \cap \Delta_n)>0,$$
and (\ref{11}) then implies that $1/2$ is an essential value for the modulus of $a$. \carre

\end{secc}

\vspace{0.2cm}

\noindent {\sc Appendix:} {\sl Proof of proposition  \ref{for}}. 

\vspace{0.2cm}

We sketch here the proof given in \cite{forrest} of proposition  \ref{for}. For the bound on $\|2q_nx\|$ note that for any strictly increasing sequence of integers $l_n$, the set of $x$ such that the sequence ${( l_n x)}_{n \in \N}$ is dense has  full Lebesgue measure. Hence we just have to show that for any infinite subset $Q$ of the set of approximation denominators of $\a$  we have that for Lebesgue almost  every $x \in [0,1]$, there exists a sequence $q_n \in Q$ such that $\lim |a(x,q_n)|=\infty$. First, it is easy to see  that the set of $x \in [0,1]$ satisfying the above condition is invariant by translation by $\a$, but $x \mapsto x+\a$ is ergodic, hence it is enough to prove that the set in question has positive measure. Next, by a simple computation we obtain for a given $k \in \N$ and any sequence $q_n$ such that $q_n \|q_n \a\| \to 0$, $2 \max \lbrace |a(x+k\a,q_n)|,|a(x,q_n)| \rbrace \geq \|2q_nx\| |a(x,k)| + C_k u_n$ where $u_n \to 0$ as $n \to \infty$ (cf. \cite[Corollary A.5]{forrest}). Hence the proof of the proposition is reduced to the following
\begin{prop}\cite[Proposition 3.13]{forrest} \label{pop} Suppose $\theta \in [0,1] \setminus \Q$ has a continued fraction representation $[a_1,a_2,\ldots]$ such that $\sum_n 1/a_n < \infty$. Then there is a $\rho >0$ such that for all $C>0$, there is a $k$ such that $\lambda \lbrace x : |a(x,k)|\geq C \rbrace \geq \rho$.
\end{prop}
To prove proposition \ref{pop}, it is convenient to define first the following function similar to the modulus of the Weyl sums
$$\psi(\a,x,k) := \left| \sum_{j=0}^{k-1} e(j^2\a/2+ jx)\right| $$
that satisfies for $0<\a$ and $x<1$
\begin{eqnarray} \label{o} \sqrt{\a} \psi(\a,x,k) = \psi(\{1/\theta \}, \{-x/\a \}, [k\a]) +O(1) \end{eqnarray}
where $\{\cdot\}$ and $[\cdot]$ denote the fractional and the integer part of a number and where the constant involved in the Landau's error notation is absolute. Equation (\ref{o})  is the only ``hard analysis'' esitmate that is needed in \cite{forrest}, but  is really crucial since it is at the center of the proof of proposition \ref{pop}. It was obtained by Hardy and Littelwood  \cite[2.128, 2.17]{HL} as a generalisation of a formula of Lindel\"of in the case of $\a$ rational and its proof is based on the calculus of residues.
 
We explain now how (\ref{o}) is used to prove proposition \ref{pop}. Given $k \in \N$, let $S\a = \{1 / \a\}$ and define $\tilde{S}(\a,x)=(S\a,\{-x/\a\})$ and write $(S^m\a,U_\a^{(m)}x) = \tilde{S}^m(\a,x)$. Let $\sigma_m(\a)=\sqrt{S^{m-1}\a} \sigma_{m-1}(\a)$ with $\sigma_0(\a)=1$, and $k(m)= [k(m-1)S^{m-1}(\a)]$ with $k(0)=k$. We have by induction from (\ref{o})
\begin{eqnarray} \label{uu} \sigma_m(\a)\psi(\a,x,k) = \psi(S^m\a,U_\a^{(m)}x,k(m))  + O(1) \end{eqnarray}
(the constant in $O(1)$  is absolute and comes from the fact that $O(\sum_{l=1}^{m} \sigma_{m-l}(S^l\a)) = O(\sum_{l=1}^m 2^{-l/2})$ since  the hypothesis $\sum 1/a_n <\infty$  implies $\lim \sup a_n \geq 2$ which in its trun implies that $\sigma_j(S^p\a) \leq C(\a) 2^{-j/2}$ for any $p$ and $j$).   

Recall the notation $b(x,k)= \sum_{j=0}^{m-1} e(jx).$
Since $\psi(\a,x,k)= |b(x,k)|+O(k^3\|\a\|)$ with an absolute constant in the error term, we have from (\ref{uu})
\begin{eqnarray} \label{uuu} \sigma_m(\a)\psi(\a,x,k) \geq |b(U_\a^{(m)}x,k(m))| - {C}({k(m)}^3\|S^m\a\|+1), \end{eqnarray}
for some absolute constant $C$.
On another hand, the condition $\sum 1/a_n < \infty$ is crucial (see \cite[Corollary 3.6]{forrest}) in checking that for all $0<\eta<1/2$ and for all $m \geq 1$
 $$\lambda\{x : \|U_\a^{(m)}x\|<\eta\} \geq \tilde{C}\eta, \quad {\rm for \ some \ absolute \ constant \ } \tilde{C}$$
which by an elementary computation implies that for any $C_0 \geq 1$
\begin{eqnarray} \lambda \{x : |b(U_\a^{(m)}x,[2\pi C_0]+1)| \geq C_0 \} &\geq& \tilde{C}/(2[2\pi C_0]+2).  \label{u} \end{eqnarray}

Fix now $C_0 \geq 3C$ where $C$ is the constant of (\ref{uuu}). Given any $C'>0$ pick $m$ sufficently large so that $C_0/(3\sigma_m(\a)) \geq C'$ and ${([2\pi C_0]+1)}^3\|S^m\a\| \leq 1$ (possible due to the arithmetical condition on $\a$). 
Let $k=k(0)$ be such that $k(m)=[2\pi C_0]+1$. We have from (\ref{uuu}) that
\begin{eqnarray*} \{ x : \psi(\a,x,k) \geq C' \}  \subset \{ x : \sigma_m(\a) \psi(\a,x,k) \geq C_0/3 \}  \subset \{ x : |b(U_\a^{(m)}x,[2\pi C_0]+1)| \geq C_0 \} \end{eqnarray*} 
and the latter set has by (\ref{u}) a measure greater than the constant  $\rho=\tilde{C}/(2[2 \pi C_0]+2). \ \ \ \ $ 

 \carre

\vspace{0.6cm}

\noindent {\bf Acknowledgments.} I am grateful to Mariusz Lema\'nczyk, Fran\c{c}ois Parreau, Jean-Paul Thouvenot and Robert Vaughan for many useful conversations,  and to Mahesh Nerurkar and Dalibor Voln\'y for communicating to me their preprint \cite{volny}. I thank the referee for his useful remarks.

\vspace{0.6cm}

\frenchspacing
\bibliographystyle{plain}

\vspace{0.4cm}

\end{document}